\documentclass[11pt]{article}
\usepackage{amsmath,amssymb,bm,color,tikz,mathrsfs,setspace}
\usetikzlibrary{matrix}
\tikzset{tab/.style={matrix of math nodes,column sep=-.35, row sep=-.35,text height=7pt,text width=7pt,align=center,inner sep=2,font=\footnotesize}}
\newcommand{\tikzmark}[2]{\tikz[overlay,remember picture,baseline] \node [anchor=base] (#1) {$#2$};}
\newcommand{\DrawLine}[3][]{%
  \begin{tikzpicture}[overlay,remember picture]
    \draw[#1] (#2.210) -- (#3.30);
  \end{tikzpicture}
}
\usepackage{hyperref}
\usepackage[margin=1.25in]{geometry}

\newcommand{\bc}{\mathbf{C}}

\newcommand{\bq}{\mathbf{Q}}

\newcommand{\bz}{\mathbf{Z}}
\newcommand{\CB}{\mathscr{B}}
\newcommand{\cc}{{\bm c}}
\newcommand{\ii}{{\bf i}}
\newcommand{\stack}[1]{\begin{smallmatrix}#1\end{smallmatrix}}
\newcommand{\TT}{\mathcal{T}}

\numberwithin{equation}{section}
\numberwithin{figure}{section}
\numberwithin{table}{section}
\newtheorem{Theorem}[equation]{Theorem}
\newtheorem{ThmDef}[equation]{Theorem/Definition}
\newtheorem{Proposition}[equation]{Proposition} 
\newtheorem{Lemma}[equation]{Lemma}

\newtheorem{Definition}[equation]{Definition}
\newtheorem{Example}[equation]{Example}
\newtheorem{Remark}[equation]{Remark}

\author{Ben Salisbury\thanks{Partially supported by CMU Early Career grant \#C62847.} \\
\small Department of Mathematics\\ 
\small Central Michigan University\\ 
\small Mount Pleasant, MI 48859, USA \\ 
\small \texttt{ben.salisbury@cmich.edu}
\and Adam Schultze\thanks{Partially supported by NSF grant DMS-1265555.} \\
\small Department of Mathematics and Statistics\\ 
\small State University of New York at Albany\\ 
\small Albany, NY 12222, USA \\
\small \texttt{aschultze@albany.edu}
\and Peter Tingley\thanks{Partially supported by NSF grant DMS-1265555.} \\
\small Department of Mathematics and Statistics\\ 
\small Loyola University of Chicago\\ 
\small Chicago, IL 60660, USA\\
\small \texttt{ptingley@luc.edu}
}
\title{\bf Combinatorial descriptions of the crystal structure on certain PBW bases (extended abstract)}

\date{\today}
\begin{document}
\maketitle
\begin{abstract}
%\paragraph{Abstract.}
Lusztig's theory of PBW bases gives a way to realize the crystal $B(\infty)$ for any simple complex Lie algebra where the underlying set consists of Kostant partitions. In fact, there are many different such realizations, one for each reduced expression for the longest element of the Weyl group. There is an algorithm to calculate the actions of the crystal operators, but it can be quite complicated. For ADE types, we give conditions on the reduced expression which ensure that the corresponding crystal operators are given by simple combinatorial bracketing rules.  We then give at least one reduced expression satisfying our conditions in every type except $E_8$, and discuss the resulting combinatorics.  Finally, we describe the relationship with more standard tableaux combinatorics in types $A$ and $D$. 
%
%\paragraph{R\'{e}sum\'{e}.}
%La th\'eorie de Lusztig de bases PBW donne un moyen de r\'ealiser le cristal $B(\infty)$ pour toute alg\`ebre de Lie complexe-simple o\`u l'ensemble sous-jacent est constitu\'e de partitions Kostant.  En fait, il existe de nombreuses r\'ealisations diff\'erentes, une
%pour chaque expression r\'eduite de l'\'el\'ement le plus long du groupe de Weyl.  Il existe un algorithme pour calculer les actions des op\'erateurs de cristal, mais elle peut \^etre assez compliqu\'e.  Pour les types ADE, nous donnons des conditions sur l'expression r\'eduite qui assurent que les op\'erateurs de cristal correspondants sont donn\'es par une simple r\`egle de d'appariement combinatoire.  Nous donnons ensuite au moins une expression r\'eduite satisfaite nos conditions pour chaque type, sauf $E_8$, et discutons la combinatoire r\'esultant. Enfin, nous d\'ecrivons les relations avec la combinatoire de tableaux standards dans les types $A$ et $D$.
\end{abstract}

\section{Introduction}

We consider the crystal $B(\infty)$ for a simple Lie algebra over $\bc$ of simply-laced type. This is a combinatorial object that contains a great deal of information about the algebra and its finite-dimensional representations. It is usually defined by a complicated algebraic construction, but it can often be realized in quite simple ways. 

Lusztig's early algebraic construction of the canonical basis of the quantum group $U_q(\mathfrak{g})$ (see \cite[Chapters 41 and 42]{L10} and references therein) can be interpreted as giving a number of parameterizations of $B(\infty)$, one for each reduced expression for the longest word $w_0$ in the Weyl group. For each of these realizations, at least one of the crystal operators is very simple, but others may be complicated. However, Lusztig explicitly describes how the realizations are related for reduced expressions that differ by a braid move. This gives a way to realize the whole crystal: Fix a reduced expression for $w_0$. Then an element of the crystal is simply a PBW monomial. To apply a crystal operator, modify the reduced expression via a sequence of braid moves until that operator is simple, then apply the operator, then modify it back. For more on this point of view, see \cite{tingley}.

This procedure is algorithmic, but can be complicated. In type $A_n$, there is a simpler realization using multisegments. As discussed in \cite{CT15}, this is precisely Lusztig's crystal structure for the reduced expression 
$
w_0= (s_1 s_2 \cdots s_n)(s_1s_2 \cdots s_{n-1}) \cdots (s_1 s_2) s_1.
$
This is not true for other reduced expressions (unless they are related to this one by trivial braid moves). 

In the current work, we give a set of conditions on a reduced expression that ensure Lusztig's crystal structure on Kostant partitions is given by a simple bracketing procedure, similar to the type $A_n$ structure on multisegments. There is at least one such reduced expression in every simply-laced type except $E_8$.

We then discuss the type $A_n$ and $D_n$ situation in some detail. There is already a nice combinatorial realization of $B(\infty)$ in these types, due to J.\ Hong and H.\ Lee \cite{HL08}, where the underlying set consists of marginally large tableaux.  We explicitly describe the unique crystal isomorphism between marginally large tableaux and Kostant partitions in these types. The isomorphism in type $A_n$ is essentially the same as the one given in \cite{CT15}, and naturally factors through the realization of $B(\infty)$ in terms of multisegments. The type $D_n$ isomorphism is new; in particular, it does not agree with the bijection used in \cite{LS14}.

The type $D_n$ situation has one interesting new feature: it is not possible to find a reduced expression that is adapted both for calculating both the ordinary crystal operators and the $*$-operators (i.e., the crystal operators twisted by Kashiwara's involution). Thus, while the ordinary crystal operators can be described with a bracketing procedure, the $*$-operators are more complicated. This may explain why the embeddings $B(\lambda) \lhook\joinrel\longrightarrow B(\infty)$ of the finite crystals into the infinity crystal seem more difficult to understand in type $D_n$, since the conditions describing which elements of $B(\infty)$ are present in a fixed $B(\lambda)$ reference the $*$-crystal structure. See \cite{HL14} for some recent work discussing these embeddings.

Reineke \cite{reineke} has also given an explicit description of the crystal operators on Lusztig's PBW basis for certain reduced expressions of $w_0$, using quiver representation theory. There the reduced expression must be adapted to an orientation of the Dynkin diagram. Our construction works nicely when the reduced expression is ``$i$-semi-adapted for all $i$'' (see \S\ref{s:sab}). These conditions do not coincide. For instance, the reduced expression we use in \S\ref{ss:etD} for type $D_n$ is not adapted to any orientation of the Dynkin diagram. 

The results here will appear with full proofs in two upcoming papers. See \cite{SST1} for the general results and \cite{SST2} for the connection to marginally large tableaux in type $D_n$.

\section{Background}

\subsection{General}

Let $\mathfrak{g}$ be a complex-simple Lie algebra of type ADE. Let $I$ be the index set of $\mathfrak{g}$, $A = (a_{ij})$ the Cartan matrix, $\{ \alpha_i \}_{i \in I}$ the positive simple roots, $\{\alpha_i^\vee\}_{i\in I}$ the simple coroots, $\Phi^+$ the set of positive roots, $P$ the weight lattice, $P^\vee$ the dual weight lattice, $W$ the Weyl group with longest element $w_0$, and $\{ s_i \}_{i \in I}$ the generating reflections.  Let $(-|-)$ denote a symmetric bilinear form on $P$ satisfying $(\alpha_i|\alpha_j) = a_{ij}$ and let $\langle-,-\rangle \colon P^\vee \times P \longrightarrow \bz$ the canonical pairing.  Let $R(w_0)$ denote the set of reduced expressions for the longest element $w_0$ of the Weyl group.  

Let $U_q(\mathfrak{g})$ be the quantized universal enveloping algebra of $\mathfrak{g}$, which is a $\bq(q)$-algebra generated by $E_i$, $F_i$, and $q^h$, for $i\in I$ and $h\in P^\vee$, subject to certain relations (see, for example, \cite{HK02}). Let $U_q^-(\mathfrak{g})$ be the subalgebra generated by the $F_i$'s. 

\subsection{Crystals}

Let $e_{i}$, $f_{i}$ be the Kashiwara operators on $U_q^-(\mathfrak{g})$ defined in \cite{K91}, $\mathcal{A} \subset \bq(q)$ be the subring of functions regular at $q=0$ and define $L(\infty)$ to be the $\mathcal{A}$-lattice spanned by
\[
S= \{  f_{i_1}  f_{i_2} \cdots  f_{i_t} \cdot 1 \in U_q^-(\mathfrak{g}) : t\ge 0, \ i_k \in I \}.
\]

\begin{ThmDef}[\cite{K91}] \hfill 
\begin{enumerate}
\item Let $\pi \colon L(\infty) \longrightarrow L(\infty)/qL(\infty)$ be the natural projection and set $B(\infty) = \pi(S)$. Then $B(\infty)$ is a $\bq$-basis of $L(\infty)/qL(\infty)$.
\item The operators $e_i$ and $f_i$ act on $L(\infty)/qL(\infty)$ for each $i\in I$. Moreover, $e_i \colon B(\infty) \longrightarrow B(\infty)\sqcup \{0\}$ and $f_i\colon B(\infty) \longrightarrow B(\infty)$ for each $i\in I$. For $b,b'\in B(\infty)$, we have $f_i b = b'$ if and only if $e_ib' = b$. 
\end{enumerate}
\end{ThmDef}

\subsection{Reduced expressions, convex orders, and PBW bases}

\begin{Definition}
A total order $\prec$ on $\Phi^+$ is called {\it convex} if, for all triples of roots $\beta, \beta',\beta''$ with $\beta'=\beta+\beta''$, we have either $\beta\prec\beta'\prec\beta''$ or $\beta''\prec\beta'\prec\beta$.
\end{Definition}

\begin{Theorem}[\cite{papi}]
There is a one-to-one correspondence between $R(w_0)$ and convex orders on $\Phi^+$.  Explicitly, if $w_0 = s_{i_1}s_{i_2}\cdots s_{i_N}$, then the corresponding convex order $\prec$ 
is defined by
\[
\beta_1 = \alpha_{i_1} \ \ \prec \ \ 
\beta_2 = s_{i_1}(\alpha_{i_2}) \ \ \prec \ \ \cdots \ \ \prec \ \ 
\beta_N = s_{i_1}s_{i_2}\cdots s_{i_{N-1}}(\alpha_{i_N}).
\]
\end{Theorem}

For $c\in \bz_{> 0}$ define
\[
F_i^{(c)} := \frac{F_i^c}{[c]!}, \  \text{where} \ \ [c]! := \prod_{j=1}^c
\frac{q^{j}-q^{-j}}{q-q^{-1}}.
\]
Given $\ii = (i_1,\dots,i_N)\in R(w_0)$ and $\cc = (c_\beta^\ii \in \bz_{\ge0}^N : \beta\in\Phi^+)$, define
\begin{equation}\label{eq:lustparam}
F_\ii^{\cc} = F_{\beta_1}^{(c_{\beta_1}^\ii)}F_{\beta_2}^{(c_{\beta_2}^\ii)}\cdots F_{\beta_N}^{(c_{\beta_N}^\ii)}
\ \ \ 
\text{ where }
\ \ \ 
F_{\beta_k}^{(c_{\beta_k}^\ii)} = T_{i_1}T_{i_2}\cdots T_{i_{k-1}}(F_{i_k}^{(c_{\beta_k}^\ii)}),
\end{equation}
and $T_i$ is the Lusztig automorphism of $U_q(\mathfrak{g})$ defined in \cite[\S 37.1.3]{L10} (where it is denoted $T''_{i,-1}$). Then the set 
\begin{equation}
\CB_\ii = \{ F_\ii ^{\cc} : \cc \in \bz_{\ge0}^N \}
\end{equation}
forms a $\bq(q)$-basis of $U_q^-(\mathfrak{g})$ called the {\it PBW basis}.  

\begin{Theorem}[\cite{saito}]
For $\ii \in R(w_0)$, 
$\mathrm{Span}_{\mathcal{A}} (\CB_\ii)  = L(\infty)$ and
$\CB_\ii+ q L(\infty)=B(\infty)$. 
\end{Theorem}

Since $\CB_\ii + qL(\infty) = B(\infty)$ for every $\ii \in R(w_0)$, there is a parametrization of the elements of $B(\infty)$ dependent on a chosen $\ii$.  Given $b \in B(\infty)$, denote the bijection which takes $b \mapsto F_\ii^\cc \bmod qL(\infty) \mapsto \cc$ by $\cc^\ii$, and call $\cc = \cc^\ii(b)$ the {\it $\ii$-Lusztig data} of $b$.  The crystal structure on $B(\infty)$ can be interpreted using these parameterizations.

\begin{Proposition}[\cite{BZ01,L10}]
Fix some $i \in I$ and let $\ii \in R(w_0)$ be such that $i_1 = i$.  Suppose $b \in B(\infty)$ and $\cc^\ii(b) = (c_1,c_2,\dots,c_N)$.  Then $\cc^\ii(f_i b) = (c_1+1,c_2,\dots,c_N)$.  Moreover, if $c_1 = 0$, then $e_i b =0$.  Otherwise, $\cc^\ii(e_i b) = (c_1-1,c_2,\dots,c_N)$.
\end{Proposition}

If $\ii$ is such that $i_1 = i$, then the operator $f_i$ is easily understood on $B(\infty)$, as indexed by $\CB_\ii$.  To understand the whole crystal structure, it is necessary to understand how the Lusztig data $\cc^\ii(b)$ are related for different $\ii$. All reduced expressions for $w_0$ are related by a sequence of braid moves, so it is enough to understand what happens to $\cc^\ii(b)$ when $\ii$ changes by a single braid move. The following is due to Lusztig. 

\begin{Lemma}[{\cite[\S 2.1]{L90}}]
\label{Lem:moves} 
Fix $\ii \in R(w_0)$ and $b \in \CB_\ii$.  Let $\{\beta_1 \prec \cdots \prec \beta_N\}$ be the order on $\Phi^+$ corresponding to $\ii$.
\begin{enumerate}
\item If $\ii'\in R(w_0)$ is such that $\ii$ and $\ii'$ differ by replacing two consecutive entries $(i_k,i_{k+1})$, such that $a_{i_k, i_{k+1}}=0$ by $(i_{k+1},i_{k})$, then $c^{\ii'}_\beta(b) = c^{\ii}_\beta(b)$ for all $\beta$ (although the order has changed).
\item If $\ii'\in R(w_0)$ is such that $\ii$ and $\ii'$ differ by replacing three consecutive $(i_k,i_{k+1},i_{k+2})$ such that $a_{i_k,i_{k+1}} = -1$ and $i_{k+2}=i_k$ by $(i_{k+1},i_k,i_{k+1})$, then
\begin{itemize}
\item $c^{\ii'}_{\beta_k}(b) = \max \{ c^\ii_{\beta_{k+1}}(b), c^\ii_{\beta_{k}}(b)+c^\ii_{\beta_{k+1}}(b)-c^\ii_{\beta_{k+2}}(b) \}$,
\item $c^{\ii'}_{\beta_{k+1}}(b) = \min \{ c^\ii_{\beta_k}(b), c^\ii_{\beta_{k+2}}(b) \}$,
\item $c^{\ii'}_{\beta_{k+2}}(b) = \max \{ c^\ii_{\beta_{k+1}}(b), c^\ii_{\beta_{k+1}}(b)+c^\ii_{\beta_{k+2}}(b)-c^\ii_{\beta_{k}}(b) \}$,
\end{itemize}
and for all other $\beta$, $c^{\ii'}_\beta(b)=c^\ii_\beta(b)$. 
\end{enumerate}
%Here, for example, $c^{\ii'}_{\beta_{k}}(b)$ means the Lusztig data for reduced expression $\ii'$, with the $k^{th}$ root in the order for $\ii$.  This could also be denoted $c^{\ii'}_{\beta'_{k+2}}(b)$, where ${\beta'_{k+2}}$ means root $k+2$ in the $\ii'$ order. 
Caution: $\beta_k$ and $\beta_{k+2}$ are actually the $(k+2)^{nd}$ and $k^{th}$ roots respectively in the order for $\ii'$.
See Example \ref{ex:pbw_compute}. 
\end{Lemma}

\begin{Lemma}
\label{Lem:braid_on_roots}
Let $\ii,\ii' \in R(w_0)$ and suppose $\{\beta_1 \prec \cdots \prec \beta_N\}$ and $\{\beta_1 \prec' \cdots \prec' \beta_N\}$ are the convex orderings on $\Phi^+$ determined by $\ii$ and $\ii'$, respectively.
\begin{enumerate}
\item The reduced expressions $\ii,\ii'$ are related by a $2$-term braid move $(i_k,i_{k+1}) \to (i_{k+1},i_k)$ if and only if $(\beta_k|\beta_{k+1})=0$.    In this case, $\beta_k \prec \beta_{k+1}$ is replaced by $\beta_{k+1} \prec' \beta_k$ after the braid move.
\item The reduced expressions $\ii,\ii'$ are related by a $3$-term braid move $(i_k,i_{k+1},i_{k+2}) \to (i_{k+1},i_k,i_{k+1})$, with $i_k = i_{k+2}$, if and only if $\{\beta_k,\beta_{k+1},\beta_{k+2}\}$ form a root system of type $\mathfrak{sl}_3$.  In this case, $\beta_k \prec \beta_{k+1} \prec \beta_{k+2}$ is replaced by $\beta_{k+2} \prec' \beta_{k+1} \prec' \beta_k$ after the braid move.
\end{enumerate} 
\end{Lemma}

For any two reduced expressions, we can understand the map $R_{\ii}^{\ii'}\colon \bz_{\ge0}^N \longrightarrow \bz_{\ge0}^N$ sending $\cc^\ii(b)$ to $\cc^{\ii'}(b)$ by finding a way to move from $\ii$ to $\ii'$ by a sequence of braid moves, and composing the maps above. 
Putting this together gives a realization of $B(\infty)$ where the underlying set is $\CB_\ii$ for some fixed $\ii$, and the $f_i$ are calculated as in the following example. 

\begin{Example}\label{ex:pbw_compute}
Let $\mathfrak{g}$ be of type $D_4$, $\ii = 123421234234$, where $\alpha_2$ is the simple root at the trivalent node of the Dynkin diagram.  The corresponding order on positive roots is
\[
1 \prec 12 \prec 123 \prec 124 \prec 1234 \prec 12234 \prec 2 \prec 24 \prec 23 \prec 234 \prec 3 \prec 4,
\]
where $1$ is identified with $\alpha_1$, $12$ with $\alpha_1 + \alpha_2$, and so on. Consider 
\begin{equation}
b= F_1^{(2)}  F_{12}^{(1)}  F_{123}^{(4)} F_{124}^{(2)} F_{1234}^{(1)} F_{12234}^{(3)} F_{2}^{(3)}  F_{24}^{(1)} F_{23}^{(2)} F_{234}^{(1)} F_{3}^{(2)} F_{4}^{(0)} \in \CB_\ii,
\end{equation}
Calculating $f_1 b$ is easy: the exponent of $F_1$ simply increases by $1$. The calculation of $f_4 b$, goes as follows:
\[
\arraycolsep=5pt
\begin{array}{rccccccccccccccc}
b &= & F_1^{(2)} & F_{12}^{(1)} & F_{123}^{(4)} & F_{124}^{(2)} & F_{1234}^{(1)} & F_{12234}^{(3)} & F_{2}^{(3)} & F_{24}^{(1)} & F_{23}^{(2)} & F_{234}^{(1)} & F_{3}^{(2)} & F_{4}^{(0)} \\[5pt]
&\simeq & F_1^{(2)} & F_{12}^{(1)} & {\color{red}F_{124}^{(2)}} & {\color{red}F_{123}^{(4)}} & F_{1234}^{(1)} & F_{12234}^{(3)} & F_{2}^{(3)} & F_{24}^{(1)} & F_{23}^{(2)} & F_{234}^{(1)} & {\color{red}F_{4}^{(0)}} & {\color{red}F_{3}^{(2)}} \\[5pt]
&\simeq & F_1^{(2)} & F_{12}^{(1)} & F_{124}^{(2)} & F_{123}^{(4)} & F_{1234}^{(1)} & F_{12234}^{(3)} & F_{2}^{(3)} & F_{24}^{(1)} & {\color{red}F_{4}^{(1)}} & {\color{red}F_{234}^{(0)}} & {\color{red}F_{23}^{(3)}} & F_{3}^{(2)} \\[5pt]
&&&&&&& \vdots \\[5pt]
 &\simeq & {\color{red}F_4^{(2)}} & {\color{red} F_{1}^{(2)}} & F_{124}^{(1)} & F_{12}^{(2)} & F_{1234}^{(1)} & F_{123}^{(4)} & F_{12234}^{(3)} & F_{24}^{(1)} & F_{2}^{(3)} & F_{234}^{(0)} & F_{23}^{(3)} & F_{3}^{(2)} \\[5pt]
f_4 b &\simeq & {\color{blue}F_4^{(3)}} & F_{1}^{(2)} & F_{124}^{(1)} & F_{12}^{(2)} & F_{1234}^{(1)} & F_{123}^{(4)} & F_{12234}^{(3)} & F_{24}^{(1)} & F_{2}^{(3)} & F_{234}^{(0)} & F_{23}^{(3)} & F_{3}^{(2)} \\[5pt]
&&&&&&& \vdots \\[5pt]
&\simeq & F_1^{(2)} & F_{12}^{(1)} & {\color{blue}F_{123}^{(3)}} & F_{124}^{(2)} & {\color{blue}F_{1234}^{(2)}} & F_{12234}^{(3)} & F_{2}^{(3)} & F_{24}^{(1)} & F_{23}^{(2)} & F_{234}^{(1)} & F_{3}^{(2)} & F_{4}^{(0)}.
\end{array}
\]
The first step performs two $2$-term braid moves and the corresponding (trivial) piecewise linear bijections. The second step is the piecewise-linear bijection for a $3$-term braid move. Specifically, by Lemma \ref{Lem:moves},
\begin{align*}
c_{23}^{\ii'}(b) &= \max\{ c_{234}^\ii(b), c_{23}^\ii(b) + c_{234}^\ii(b) - c_4^\ii(b) \} = \max\{ 1,2+1-0\} = 3 ,\\
c_{234}^{\ii'}(b) &= \min\{ c_{23}^\ii(b), c_{4}^\ii(b) \} = \min\{ 2,0\} = 0,\\
c_{4}^{\ii'}(b) &=  \max\{ c_{234}^\ii(b), c_{234}^\ii(b) + c_4^\ii(b) - c_{23}^\ii(b) \} = \max\{ 1, 1+0-2\} = 1.
\end{align*}
This, along with the reordering of the roots by Lemma \ref{Lem:braid_on_roots}, gives the third line above.  Continue making braid moves and applying the corresponding piecewise-linear bijections to get a PBW monomial with $F_4$ as the leftmost factor. (To do this, recall that, by Lemma \ref{Lem:braid_on_roots}, two roots can be interchanged with a $2$-term braid move exactly if they are perpendicular, and that a three-term braid move applies to three consecutive roots $\beta, \beta',\beta''$ if and only if $\beta'=\beta+\beta''$.) Then increase that exponent by $1$. Then do braid moves and the corresponding piecewise-linear bijections to get back to the original order (not shown). The result is $f_4b$, expressed as an element of $\CB_\ii$. 
\end{Example}

\section{Semi-adapted words and bracketing} \label{s:sab}

In general, calculating $f_i$ as in Example \ref{ex:pbw_compute} can be computationally involved. However, for some words, the application of $f_i$ can be calculated by a simple bracketing procedure. We now discuss those words. 

\subsection{Semi-adapted words}

\begin{Definition} \label{def:sa}
Fix a reduced expression ${\bf i}$ for $w_0$, and $i \in I$. We say that $\ii$ is {\it adapted} for $i$ if $i_1=i$. We say that $\ii$ is {\it semi-adapted} for $i$ if one can perform a sequence of braid moves to $\ii$ to get to a word $\ii'$ with $i'_1=i$, and each of these is either
\begin{itemize}
\item a $2$-term braid move, or
\item a $3$-term braid move such that the corresponding roots before the move are $(\beta, \beta+\alpha_i, \alpha_i)$, in that order, for some $\beta$. 
\end{itemize}

\end{Definition} 

\begin{Definition} \label{def:em}
Fix an $i$-semi-adapted reduced expression $\ii$. Let $\eta_1, \ldots, \eta_k$ be the roots prior to $\alpha_i$ in the corresponding order such that $(\alpha_i|\eta_j) <0$, in the order they appear left to right. Let $\nu_j=\eta_j+\alpha_i$, which is a root. 
\end{Definition}

\subsection{Bracketing rules}

Recall that a Kostant partition is a formal $\bz_{\ge0}$-linear combination of positive roots, which we denote by $(c_\beta)_{\beta\in\Phi^+}$.  Therefore, an $\ii$-Lusztig datum $\cc^\ii(b)$ may be identified with a Kostant partition.

\begin{Definition}\label{def:bso}
For any Kostant partition $\cc = (c_\beta)_{\beta\in\Phi^+}$, define $S_i^\ii(\cc)$ to be the string of brackets
\begin{center}
\begin{tikzpicture}
\node at (0,0) {$) \cdots )$};
\node at (1,0) {$( \cdots ($};
\node at (2,0) {$\cdots$};
\node at (3,0) {$) \cdots )$};
\node at (4,0) {$( \cdots ($};
\node at (5.3,0) {$) \cdots ) \; .$};

\node at (0,0.5) {$c_{\nu_1}$};
\node at (1,0.5) {$c_{\eta_1}$};
\node at (3,0.5) {$c_{\nu_k}$};
\node at (4,0.5) {$c_{\eta_k}$};
\node at (5.3,0.5) {$c_{\alpha_i}$};
\end{tikzpicture}
\end{center}
We sequentially pair adjacent brackets $()$ until the remaining brackets are a subsequence of the form $))\cdots)((\cdots($.  A bracket in $S_i^\ii(\cc)$ is called {\it uncanceled} if it is not paired in this procedure.
\end{Definition}

\begin{Theorem} \label{th:mb}
Fix $i \in I$ and an $i$-semi-adapted word $\ii$. Let $b \in \CB_\ii$, and let $\cc = (c_\beta^\ii(b))$ be the corresponding Kostant partition. 
\begin{itemize}
\item If the leftmost uncanceled `$($' in $S_i^\ii(\cc)$ corresponds to $\eta_j$ then 
\[
c_{\nu_j}^\ii(f_i b)= c_{\nu_j}^\ii(b)+1,\qquad
 c_{\eta_j}^\ii(f_i b)= c_{\eta_j}^\ii(b)-1, \quad \text{and } \quad
 c_\beta^\ii(f_i b)=c_\beta^\ii(b) \ \text { for all other $\beta$}.
\]
If there is no uncanceled `$($', then $c_{\alpha_i}^\ii(b)$ increases by $1$ and all other $c_\beta^\ii(b)$ are unchanged. 
\item If the rightmost uncanceled `$)$' in $S_i^\ii(\cc)$ corresponds to $\nu_j$ then 
\[
c_{\nu_j}^\ii(e_i b)= c_{\nu_j}^\ii(b)-1, \qquad
c_{\eta_j}^\ii(e_i b)= c_{\eta_j}^\ii(b)+1, \quad \text{and} \quad
c_\beta^\ii(e_i b)=c_\beta^\ii(b) \text{ for all other $\beta$.}
\]
If there is no uncanceled `$)$' in $S_i^\ii(\cc)$, then $e_ib = 0$.
\end{itemize}
\end{Theorem}
Thus, if $\ii \in R(w_0)$ is $i$-semi-adapted for all $i$, then the crystal structure on $\CB_\ii$ can be described completely using bracketing rules. The proof of Theorem \ref{th:mb} will appear in \cite{SST1}. It relies on the following lemma.

\begin{Lemma}
In any sequence of braid moves allowed by Definition \ref{def:sa}, the $3$-term braid moves happen in the order $(\eta_k,\nu_k, \alpha_i)$, then $(\eta_{k-1},\nu_{k-1}, \alpha_i)$, and so on. 
\end{Lemma}

\subsection{Existence of semi-adapted words} 

In \cite{Lit98}, Littelmann details specific reduced expressions for the longest element of the Weyl group for which the conditions describing the associated string cone can be expressed.  These reduced expressions are called ``nice decompositions'' and correspond to ``good enumerations'' of the underlying Dynkin diagram.  At least one such nice decomposition is given for each finite type except $E_8$ and $F_4$.  As will be shown in \cite{SST1}, the reverse of these nice decompositions are $i$-semi-adapted for all $i$.  In particular, the words $\ii^A$ and $\ii^D$ described in the subsequent sections are precisely the reverse of Littelmann's nice decompositions in types $A$ and $D$, respectively.

\section{Explicit descriptions}

\subsection{Type $A$} \label{ss:bra}

The Lie algebra of type $A_n$ has Dynkin diagram
$$
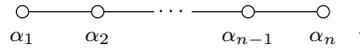
\begin{tikzpicture}[baseline=-5,scale=1,font=\scriptsize]
\foreach \x in {1,2,4,5}
{\node[circle,draw,scale=.45] (\x) at (\x,0) {};}
\node[label={below:$\alpha_1$}] at (1,0) {}; 
\node[label={below:$\alpha_2$}] at (2,0) {}; 
\node[label={below:$\alpha_{n-1}$}] at (4,0) {}; 
\node[label={below:$\alpha_{n}$}] at (5,0) {}; 
\node[label={below:.}] at (5.5,0) {}; 
\node at (3,0) {$\cdots$};
\draw[-] (1) -- (2);
\draw[-] (2) -- (2.75,0);
\draw[-] (3.25,0) -- (4);
\draw[-] (4) -- (5);
\end{tikzpicture}
$$
Write the positive roots as $\{\alpha_{i,j}: 1 \le i \le j \le n\}$, where $\alpha_{i,j} = \alpha_i + \alpha_{i+1} + \cdots + \alpha_j$.
The word $\ii^A$ corresponding to the reduced expression $w_0 = (s_1s_2\cdots s_n)(s_1s_2\cdots s_{n-1}) \cdots (s_1s_2)s_1$ is $i$-semi-adapted for all $i$.  Therefore, the crystal structure on $\CB_{\ii^A}$ is given by  bracketing rules.  

In \cite{CT15}, they show that the crystal structure on $\CB_{\ii^A}$ essentially gives the realization of $B(\infty)$ using the multisegments from \cite{JL09,LTV99}. That is, they show that the map which takes each root $\alpha_{i,j}$ to the segment $[i,j]$ is a crystal isomorphism. This isomorphism can also be understood using tableaux, which we discuss in more detail in the next section.

To understand why this map is a crystal isomorphism, consider the corresponding order on positive roots. In this case, $\alpha_{i,j} \prec \alpha_{i',j'}$ if and only if $i < i'$ or $i = i'$ and $j< j'$.  Then, given $i\in I$ and a Kostant partition $\cc = (c_\beta)_{\beta\in\Phi^+}$, the string of brackets $S_i^{\ii^A}(\cc)$ is 
\[
\underbrace{)\cdots)}_{c_{\alpha_{1,i}}}\ 
\underbrace{(\cdots(}_{c_{\alpha_{1,i-1}}}\
\underbrace{)\cdots)}_{c_{\alpha_{2,i}}}\ 
\underbrace{(\cdots(}_{c_{\alpha_{2,i-1}}}\ \ \
\cdots\ \ \
\underbrace{)\cdots)}_{c_{\alpha_{i-i,i}}}\ 
\underbrace{(\cdots(}_{c_{\alpha_{i-i,i-1}}}\
\underbrace{)\cdots)}_{c_{\alpha_{i,i}}}.
\]
The order on roots agrees with the order on segments used to define the crystal operators in \cite{JL09}.

\begin{Example}\label{ex:pbw_A}
Consider type $A_3$ with $\ii = 123121$ and $i=2$.  The corresponding order on positive roots is $1 \prec \stack{2\\1} \prec \stack{3\\2\\1} \prec 2 \prec \stack{3\\2} \prec 3$, where we identify $1$ with $\alpha_1$, $\stack{2\\1}$ with $\alpha_1+\alpha_2$, and so on.  If $b \in B(\infty)$ is such that $c^\ii(b) = (2,3,1,3,3,2)$, then the corresponding Kostant partition is 
\[
\begin{array}{cccccccccccccl}
1 & 1 & \stack{2\\1} & \stack{2\\1} & \stack{2\\1} & \stack{3\\2\\1} & 2 & 2 & 2 & \stack{3\\2} & \stack{3\\2} & \stack{3\\2}& 3 & 3.
\end{array}
\]
Placing the parts/roots in the order prescribed by Definition \ref{def:em}, we get 
\[
\arraycolsep=2pt
\begin{array}{ccccccccccccclll}
\stack{2\\1} & \stack{2\\1} & \stack{2\\1} & 1 & 1 & 2 & 2 & 2 \\[5pt]
) & ) & ) & \tikzmark{leftA1}{{\color{green}(}} & {\color{green}(} & {\color{green})}  & \tikzmark{rightA1}{{\color{green})}} & )  
\end{array}
\ \ \ \ \raisebox{2pt}{$\xrightarrow{\ \ f_2\ \ }$}\ \ \ \ 
\begin{array}{cccccccccccccl}
\stack{2\\1} & \stack{2\\1} & \stack{2\\1} & 1 & 1 & 2 & 2 & 2 & {\color{blue}2}. \\[5pt]
\phantom{(}
\end{array}
\]
\DrawLine[green, thick, opacity=0.5]{leftA1}{rightA1}

\noindent Hence, applying $f_2$ to the given partition above yields
$
\arraycolsep=2pt
\begin{array}{ccccccccccccccl}
1 & 1 & \stack{2\\1} & \stack{2\\1} & \stack{2\\1} & \stack{3\\2\\1} & 2 & 2 & 2 &  {\color{blue}2} & \stack{3\\2} & \stack{3\\2} & \stack{3\\2}& 3 & 3.
\end{array}
$
\end{Example}

However, in type $A$ there are many other good enumerations, and hence many other reduced expressions where the crystal structure on the corresponding PBW monomials is given by a bracketing rule. For example, in type $A_4$, 
$$
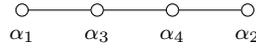
\begin{tikzpicture}[baseline=-3,scale=1,font=\scriptsize]
\foreach \x in {1,2,3,4}
{\node[circle,draw,scale=.45] (\x) at (\x,0) {};}
\node[label={below:$\alpha_1$}] at (1,0) {}; 
\node[label={below:$\alpha_3$}] at (2,0) {}; 
\node[label={below:$\alpha_4$}] at (3,0) {}; 
\node[label={below:$\alpha_2$}] at (4,0) {}; 
\draw[-] (1) -- (2);
\draw[-] (2) -- (3);
\draw[-] (3) -- (4);
\end{tikzpicture}
$$
is a good enumeration. It could be very interesting to understand the corresponding combinatorics.

\subsection{Type $D$} \label{ss:etD}

The Lie algebra of type $D_n$ has Dynkin diagram
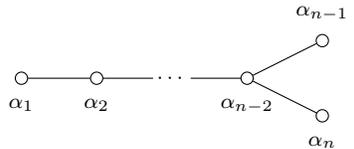
\[
\begin{tikzpicture}[baseline,scale=1,font=\scriptsize]
\foreach \x in {1,2,4}
{\node[circle,draw,scale=.45] (\x) at (\x,0) {};}
\node[label={below:$\alpha_1$}] at (1,0) {}; 
\node[label={below:$\alpha_2$}] at (2,0) {}; 
\node[label={below:$\alpha_{n-2}$}] at (4,0) {}; 
\node[label={above:$\alpha_{n-1}$}] at (5,.5) {}; 
\node[label={below:$\alpha_{n}$}] at (5,-.5) {};
\node at (3,0) {$\cdots$};
\node[circle,draw,scale=.45] (5) at (5,.5) {};
\node[circle,draw,scale=.45] (6) at (5,-.5) {};
\draw[-] (1) -- (2);
\draw[-] (2) -- (2.75,0);
\draw[-] (3.25,0) -- (4);
\draw[-] (4) -- (5);
\draw[-] (4) -- (6);
\node at (5.5,-0.85) {.};
\end{tikzpicture}
\]
The list of positive roots in type $D_n$ is given in Table \ref{posroots}.

\begin{table}[t]
\doublespacing
\[
\begin{array}{|cl|}\hline
\beta_{i,k}= \alpha_i + \cdots + \alpha_{k}, & 1\le i\le k \le n-1 \\
\gamma_{i,k}=\alpha_i + \cdots + \alpha_{n-2}+ \alpha_{n} + \alpha_{n-1}+ \cdots + \alpha_k, & 1\le i < k \le n\\\hline
\end{array}
\]
\caption{Positive roots of type $D_n$.}\label{posroots}
\end{table}

\begin{Lemma}
Define $\ii^D$ to be the word associated to the reduced expression
\[
w_0 = (s_1s_2\cdots s_{n-1}s_ns_{n-2}\cdots s_2s_1)(s_2\cdots s_{n-1}s_ns_{n-2}\cdots s_2) \cdots (s_{n-2}s_{n-1}s_ns_{n-2})s_{n-1}s_n.
\]  
Then $\ii^D$ is semi-adapted for all $i$.  
\end{Lemma}

The order of the positive roots corresponding to the subword $(i,i+1,\dots,n,n-2,\dots,i)$ of $\ii^D$, for $1\le i \le n-2$, as
\[
\begin{cases}
\beta_{i,i} \prec \beta_{i,i+1} \prec \dots \prec \beta_{i,n-2} \prec \beta_{i,n-1} \prec \gamma_{i,n} \prec \gamma_{i,n-1} \prec \cdots \prec \gamma_{i,i+1} & \text{ if } i \equiv 1 \bmod 2, \\
\beta_{i,i} \prec \beta_{i,i+1} \prec \dots \prec \beta_{i,n-2} \prec \gamma_{i,n} \prec \beta_{i,n-1} \prec \gamma_{i,n-1} \prec \cdots \prec \gamma_{i,i+1} & \text{ if } i\equiv 0 \bmod 2.
\end{cases}
\]
The ordering on the roots corresponding to the suffix $(n-1,n)$ of $\ii^D$ is 
\[
\begin{cases}
\beta_{n-1,n-1} \prec \gamma_{n-1,n} & \text{ if } n \equiv 0 \bmod 2, \\
\gamma_{n-1,n} \prec \beta_{n-1,n-1} & \text{ if } n \equiv 1 \bmod 2.
\end{cases}
\]
It follows that, for a Kostant partition $\cc = (c_\beta)_{\beta\in \Phi^+}$, the string of brackets $S_i^{\ii^D}(\cc)$ needed to compute $f_i$ is obtained by canceling brackets in
\[
\begin{cases}
\ \underbrace{)\cdots)}_{c_{\beta_{1,i}}}\ 
\underbrace{(\cdots(}_{c_{\beta_{1,i-1}}}\
\underbrace{)\cdots)}_{c_{\gamma_{1,i}}}\ 
\underbrace{(\cdots(}_{c_{\gamma_{1,i+1}}}\ 
\cdots\ 
\underbrace{)\cdots)}_{c_{\beta_{i-1,i}}}\ 
\underbrace{(\cdots(}_{c_{\beta_{i-1,i-1}}}\
\underbrace{)\cdots)}_{c_{\gamma_{i-1,i}}}\ 
\underbrace{(\cdots(}_{c_{\gamma_{i-1,i+1}}}\
\underbrace{)\cdots)}_{c_{\beta_{i,i}}}, 
& \text{ if } i \neq n, \\
\ \underbrace{)\cdots)}_{c_{\gamma_{1,n}}}\ 
\underbrace{(\cdots(}_{c_{\beta_{1,n-2}}}\
\underbrace{)\cdots)}_{c_{\gamma_{1,n-1}}}\ 
\underbrace{(\cdots(}_{c_{\beta_{1,n-1}}}\ 
\cdots\ 
\underbrace{)\cdots)}_{c_{\gamma_{n-2,n}}}\ 
\underbrace{(\cdots(}_{c_{\beta_{n-2,n-2}}}\
\underbrace{)\cdots)}_{c_{\gamma_{n-2,n-1}}}\ 
\underbrace{(\cdots(}_{c_{\beta_{n-2,n-1}}}\
\underbrace{)\cdots)}_{c_{\gamma_{n-1,n}}},
& \text{ if } i = n.
\end{cases}
\]

\begin{Example}
Consider the setup from Example \ref{ex:pbw_compute}. The corresponding Kostant partition is
\[
\arraycolsep=2pt
\begin{array}{cccccccccccccccccccccl}
1 & 1 & \stack{2\\1} & \stack{3\\2\\1} & \stack{3\\2\\1} & \stack{3\\2\\1} & \stack{3\\2\\1} & \stack{4\\2\\1} & \stack{4\\2\\1} & \stack{3\,4\\2\\1} & \stack{2\\3\,4\\2\\1} & \stack{2\\3\,4\\2\\1} & \stack{2\\3\,4\\2\\1} & 2 & 2 & 2 & \stack{4\\2} & \stack{3\\2} & \stack{3\\2} & \stack{3\,4\\2} & 3 & 3,
\end{array}
\]
where, for example, $\stack{3\,4\\2\\1}$ is the root denoted by $1234$ in Example \ref{ex:pbw_compute}; we use this new format for reasons which are explained in detail in \cite{SST1}.
Arranging the roots and making a string of brackets to calculate $f_4$ as in Definition \ref{def:bso},
\[
\arraycolsep=1pt
\begin{array}{cccccccccccccccccccclll}
\stack{4\\2\\1} & \stack{4\\2\\1} & \stack{2\\1} & \stack{3\,4\\2\\1} & \stack{3\\2\\1} & \stack{3\\2\\1} & \stack{3\\2\\1} & \stack{3\\2\\1} & \stack{4\\2} & 2 & 2 & 2 & \stack{3\,4\\2} & \stack{3\\2} & \stack{3\\2} \\[5pt]
) & ) & \tikzmark{leftD1}{\color{green}(} & \tikzmark{rightD1}{\color{green})} & {\color{blue}\boldsymbol{(}} & ( & ( & \tikzmark{leftD2}{\color{green}(} & \tikzmark{rightD2}{\color{green})} & ( & ( & \tikzmark{leftD3}{\color{green}(} & \tikzmark{rightD3}{\color{green})} & ( & ( & 
\end{array}
\ \ \ \ \raisebox{2pt}{$\xrightarrow{\ \ f_4\ \ }$}\ \ \ \ 
\begin{array}{ccccccccccccccccccccclll}
\stack{4\\2\\1} & \stack{4\\2\\1} & \stack{2\\1} & \stack{3\,4\\2\\1} & {\color{blue}\stack{3\,4\\2\\1}} & \stack{3\\2\\1} & \stack{3\\2\\1} & \stack{3\\2\\1}& \stack{4\\2} & 2 & 2 & 2 & \stack{3\,4\\2} & \stack{3\\2} & \stack{3\\2}, \\[5pt]
\phantom{(} &
\end{array}
\]
\DrawLine[green, thick, opacity=0.5]{leftD1}{rightD1}
\DrawLine[green, thick, opacity=0.5]{leftD2}{rightD2}
\DrawLine[green, thick, opacity=0.5]{leftD3}{rightD3}

\noindent which agrees with our previous calculation.
\end{Example}

\subsection{Types $E_6$, $E_7$, and $E_8$}

For types $E_6$ and $E_7$ Littelmann \cite[\S 8]{Lit98} found nice decompositions, so our results show that the crystal operators on the PBW monomials for the reverse words are given by bracketing. We did not find it enlightening to work out the details. 

In type $E_8$, Littelmann claims there is no good enumeration and hence no nice decomposition.  This does not immediately imply there is no reduced expression that is $i$-semi-adapted for all $i$, but we conjecture that no such expression exists.  Using the fact that $E_7$ is a sub-root system, it is certainly possible to find a reduced expression which is $i$-semi-adapted for all but one $i$, so all but one of the crystal operators are given by bracketing rules.  This may still be computationally useful. 

\subsection{Non-simply-laced types}

We expect similar results to hold in types $B_n$ and $C_n$ using Littelmann's nice decompositions. The main difficulty is that the rank two crystals are more complicated.    
We conjecture that in type $F_4$ there is no word that is $i$-semi-adapted for all $i$, since Littelmann finds no nice decomposition in that case.

\section{Relation to tableaux combinatorics}

\subsection{Type $A$}

Following \cite{HL08}, a marginally large tableaux of type $A_n$ is a semistandard Young tableaux $T$ on the alphabet $\{1,\dots,n+1\}$ with $n$ rows such that the number of boxes of content $i$ in the $i$th row (from the top, using the English convention) is exactly 1 more than the total number of boxes in the $(i+1)$st row. The set  $\TT(\infty)$ has a natural crystal structure. 

Given a marginally large tableau $T$, define a Kostant partition $\Theta(T)$ by setting each $c_{\alpha_i + \cdots + \alpha_j}$ to be the number of boxes with content $j+1$ on row $i$. 
It follows immediately from the results in \cite{CT15} that $\Theta$ is a crystal isomorphism between $\TT(\infty)$ and $\CB_{\ii^A}$.

\begin{Example}
Consider the setup from Example \ref{ex:pbw_A}.  Then 
\[
\Theta^{-1}(b) = 
\begin{tikzpicture}[baseline]
\matrix [tab] 
 {
	\node[draw,fill=gray!30]{1}; &
	\node[draw,fill=gray!30]{1}; &
	\node[draw,fill=gray!30]{1}; &
	\node[draw,fill=gray!30]{1}; &
	\node[draw,fill=gray!30]{1}; &
	\node[draw,fill=gray!30]{1}; &
	\node[draw,fill=gray!30]{1}; &
	\node[draw,fill=gray!30]{1}; &
	\node[draw,fill=gray!30]{1}; &
	\node[draw,fill=gray!30]{1}; &
	\node[draw,fill=gray!30]{1}; &
	\node[draw]{2}; & 
	\node[draw]{2}; &
	\node[draw]{3}; &
	\node[draw]{3}; & 
	\node[draw]{3}; & 
	\node[draw]{4};\\
	\node[draw,fill=gray!30]{2}; &
	\node[draw,fill=gray!30]{2}; &
	\node[draw,fill=gray!30]{2}; &
	\node[draw,fill=gray!30]{2}; &
	\node[draw]{3}; & 
	\node[draw]{3}; &
	\node[draw]{3}; &
	\node[draw]{4}; &
	\node[draw]{4}; &
	\node[draw]{4}; \\
	\node[draw,fill=gray!30]{3}; &
	\node[draw]{4}; &  
	\node[draw]{4}; \\
 };
\end{tikzpicture}
\ .
\]
Computing $f_3$ on this tableaux gives
\[
f_2 \Theta^{-1}(b) = 
\begin{tikzpicture}[baseline]
\matrix [tab] 
 {
	\node[draw,fill=gray!30]{1}; &
	\node[draw,fill=gray!30]{1}; &
	\node[draw,fill=gray!30]{1}; &
	\node[draw,fill=gray!30]{1}; &
	\node[draw,fill=blue]{\color{white}\bf1};&
	\node[draw,fill=gray!30]{1}; &
	\node[draw,fill=gray!30]{1}; &
	\node[draw,fill=gray!30]{1}; &
	\node[draw,fill=gray!30]{1}; &
	\node[draw,fill=gray!30]{1}; &
	\node[draw,fill=gray!30]{1}; &
	\node[draw,fill=gray!30]{1}; &
	\node[draw]{2}; & 
	\node[draw]{2}; &
	\node[draw]{3}; &
	\node[draw]{3}; & 
	\node[draw]{3}; & 
	\node[draw]{4};\\
	\node[draw,fill=gray!30]{2}; &
	\node[draw,fill=gray!30]{2}; &
	\node[draw,fill=gray!30]{2}; &
	\node[draw,fill=gray!30]{2}; &
	\node[draw,fill=blue]{\color{white}\bf3};&
	\node[draw]{3}; & 
	\node[draw]{3}; &
	\node[draw]{3}; &
	\node[draw]{4}; &
	\node[draw]{4}; &
	\node[draw]{4}; \\
	\node[draw,fill=gray!30]{3}; &
	\node[draw]{4}; &  
	\node[draw]{4}; \\
 };
\end{tikzpicture}\ .
\]
One can verify that this agrees with the calculation in Example \ref{ex:pbw_A}.
\end{Example}

\subsection{Type $D$}

As in \cite{HL08}, let $\TT(\infty)$ be the set of marginally large semistandard tableaux on the alphabet
\[
J(D_n) := \left\{ 1 \prec \cdots \prec n-1 \prec \begin{array}{c} n \\ \overline n \end{array} \prec \overline{n-1} \prec \cdots \prec \overline 1\right\},
\]
with $n-1$ rows such that the contents of the $i$th row are less than or equal to $\overline\imath$ (for each $i=1,\dots,n-1$) and entries $n$ and $\overline n$ do not appear in the same row.

\begin{figure}[h]
\[
\begin{tikzpicture}[baseline=-4,xscale=1.25,font=\footnotesize]
 \node (1) at (0,0) {$\boxed1$};
 \node (d1) at (1.5,0) {$\cdots$};
 \node (n-1) at (3,0) {$\boxed{n-1}$};
 \node (n) at (4.5,.75) {$\boxed{n}$};
 \node (bn) at (4.5,-.75) {$\boxed{\overline{n}}$};
 \node (bn-1) at (6,0) {$\boxed{\overline{n-1}}$};
 \node (d2) at (7.5,0) {$\cdots$};
 \node (b1) at (9,0) {$\boxed{\overline{1}}$};
 \draw[->] (1) to node[above]{\tiny$1$} (d1);
 \draw[->] (d1) to node[above]{\tiny$n-2$} (n-1);
 \draw[->] (n-1) to node[above,sloped]{\tiny$n-1$} (n);
 \draw[->] (n-1) to node[below,sloped]{\tiny$n$} (bn);
 \draw[->] (n) to node[above,sloped]{\tiny$n$} (bn-1);
 \draw[->] (bn) to node[below,sloped]{\tiny$n-1$} (bn-1);
 \draw[->] (bn-1) to node[above]{\tiny$n-2$} (d2);
 \draw[->] (d2) to node[above]{\tiny$1$} (b1);
\end{tikzpicture}
\]
\caption{The fundamental crystal of type $D_n$.}\label{fundD}
\end{figure}
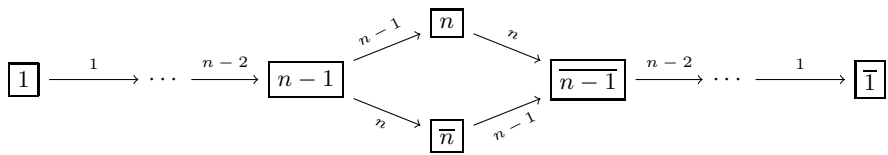

As in type $A$, $\TT(\infty)$ has a crystal structure. 
To compute the action of the Kashiwara operators $e_i$ and $f_i$ on $T \in \TT(\infty)$, 
form a sequence of brackets by reading boxes in the tableau in rows from right to left, starting with the top row. Add `)' under any letter for which there is an $i$-colored arrow entering the corresponding box in Figure \ref{fundD}, and add `(' under any letter for which there is an $i$-colored arrow leaving the corresponding box in Figure \ref{fundD}. Sequentially cancel all ()-pairs to obtain a sequence of the form $)\cdots)(\cdots($.  The brackets that remain are called {\it uncanceled}.

\begin{Definition}\label{def:T-ops}
Let $T\in \TT(\infty)$ and $i\in I$.
\begin{enumerate}
\item Let $x$ be the letter in $T$ corresponding to the right-most uncanceled `$)$.'  Then $e_iT$ is the tableau obtained from $T$ by replacing the box containing $x$ by the box containing the letter at the other end of the $i$-arrow from $x$ in Figure \ref{fundD}.  If the result is not marginally large, then delete exactly one column containing the elements $1,\dots,i$ so that the result is marginally large.  If there is no uncanceled `$)$,' then $e_iT =0$.
\item Let $y$ be the letter in $T$ corresponding to the left-most uncanceled `$($.'  Then $f_iT$ is the tableau obtained from $T$ by replacing the box containing $y$ by the box containing the letter at the other end of the $i$-arrow from $y$ in Figure \ref{fundD}.  If the result is not marginally large, then insert exactly one column containing the elements $1,\dots,i$ so that the result is marginally large.
\end{enumerate}
\end{Definition}

Definition \ref{def:T-ops} uses the middle-Eastern reading, as defined in \cite{HK02}, whereas in \cite{HL08} they use the far-Eastern reading.  However, the resulting crystal graphs are identical.

\begin{Proposition}
The crystal graphs obtained from $\TT(\infty)$ using the far-Eastern reading and the middle-Eastern reading, respectively, are identical.
\end{Proposition}

\begin{Remark}
In contrast to type $A_n$, the crystal structure on irreducible highest weight crystals of type $D_n$ modeled by Kashiwara-Nakashima tableaux \cite{KN94} using these two readings are not in fact the same. This is only a property of marginally large tableaux. 
\end{Remark}

\begin {Definition} 
Define a map $\Psi$ from $\TT(\infty)$ to Kostant partitions as follows.  
For a tableaux $T \in \TT(\infty)$, let $R_1,\dots,R_{n-1}$ denote the rows of $T$ starting at the top.  Set $\Psi(T) = (c_\beta)_{\beta\in\Phi^+}$, where $(c_\beta)$ is obtained from the trivial data $(0,0,\dots,0)$ in the following way: 
\begin{enumerate}
\item if $j \neq n-1$, each $\overline{\jmath}$ in $R_j$ increases both $c_{\beta_{j,j}}$ and $c_{\gamma_{j,j+1}}$ by $1$;

\item if $j = n-1$, each $\overline{\jmath}$ in $R_j$ increases both $c_{\beta_{n-1,n-1}}$ and $c_{\gamma_{n-1,n}}$ by $1$;

\item for each pair $k, \overline k$ in $R_j$, where $k\neq n-1$, increase both $c_{\beta_{j,k}}$ and $c_{\gamma_{j,k+1}}$ by $1$;

\item for each pair $n-1, \overline{n-1}$ in $R_j$, increase both $c_{\beta_{j,n-1}}$ and $c_{\gamma_{j,n}}$ by $1$;

\item each remaining $k\in\{j,j+1,\dots,n\}$ in $R_j$ increases $c_{\beta_{j,k-1}}$ by $1$;

\item each remaining $\overline{k}\in \{\overline{n},\overline{n-1},\dots,\overline{\jmath+1} \}$ in $R_j$ increases $c_{\gamma_{j,k}}$ by $1$.

\end{enumerate}
\end{Definition}

\begin{Example}
Let $n=4$ and 
\[
T = \begin{tikzpicture}[baseline]
\matrix [tab] 
 {
	\node[draw,fill=gray!30]{1}; &
	\node[draw,fill=gray!30]{1}; &
	\node[draw,fill=gray!30]{1}; &
	\node[draw,fill=gray!30]{1}; &
	\node[draw,fill=gray!30]{1}; &
	\node[draw,fill=gray!30]{1}; &
	\node[draw,fill=gray!30]{1}; &
	\node[draw,fill=gray!30]{1}; &
	\node[draw,fill=gray!30]{1}; &
	\node[draw]{2}; & 
	\node[draw]{2}; &
	\node[draw]{\overline 3}; &
	\node[draw]{\overline 1}; & 
	\node[draw]{\overline 1}; & 
	\node[draw]{\overline 1};\\
	\node[draw,fill=gray!30]{2}; &
	\node[draw,fill=gray!30]{2}; &
	\node[draw,fill=gray!30]{2}; &
	\node[draw,fill=gray!30]{2}; &
	\node[draw]{3}; & 
	\node[draw]{\overline 4}; &
	\node[draw]{\overline 3}; &
	\node[draw]{\overline 3}; \\
	\node[draw,fill=gray!30]{3}; &
	\node[draw]{\overline 4}; &  
	\node[draw]{\overline 3}; \\
 };
\end{tikzpicture}\ .
\]
The sequence of $\overline1$'s in the first row increases both $c_{\beta_{1,1}}$ and $c_{\gamma_{1,2}}$ by $3$, the $\overline 3$ in the first row increases $c_{\gamma_{1,3}}$ by $1$ and the pair of $2$'s in the first row increases $c_{\beta_{1,1}}$ by $2$.  In the second row, the pair $(3,\overline3)$ increases both $c_{\beta_{2,3}}$ and $c_{\gamma_{2,4}}$ by $1$, the $\overline4$ increases $c_{\gamma_{2,3}}$ by $1$, and the $3$ increases $c_{\gamma_{2,4}}$ by $1$.  Finally, the $\overline3$ in the third row increases both $c_{\beta_{3,3}}$ and $c_{\gamma_{3,4}}$ by $1$ and the $\overline4$ increases $c_{\gamma_{3,4}}$ by $1$.  In summary, 
\[
\arraycolsep=2pt
\begin{array}{rccccccccccccccccl}
\Psi(T) &=& 1 & 1 & 1 & 1 & 1 & \stack{3\,4\\2\\1} & \stack{2\\3\,4\\2\\1} & \stack{2\\3\,4\\2\\1} & \stack{2\\3\,4\\2\\1} & \stack{4\\2} & \stack{4\\2} & \stack{3\\2} & \stack{3\,4\\2} & 3 & 4 & 4.
\end{array}
\]
\end{Example}

Recall that Kostant partitions correspond to elements of $\CB_{\ii^D}$. The following is the main result of \cite{SST2}:

\begin{Theorem} 
The map $\Psi\colon \TT(\infty) \longrightarrow \CB_{\ii^D}$ defined above is a crystal isomorphism.
\end{Theorem}

\bibliography{KP-crystal}{}
\bibliographystyle{amsplain}
\end{document}